\documentclass{article}
\usepackage[a4paper,top=3cm,bottom=3.5cm,left=3.5cm,right=3.5cm,marginparwidth=1.75cm]{geometry}
\usepackage[hidelinks]{hyperref}
\usepackage{amsmath, amssymb, amsthm, bbm, comment, graphicx, cleveref, caption, siunitx, color, tikz, bigints, mathtools}
\usetikzlibrary{patterns}
\usetikzlibrary{decorations.pathmorphing}

\newcommand{\ddp}[2]{\frac{\partial#1}{\partial#2}}
\newcommand{\dd}[2]{\frac{\mathrm{d}#1}{\mathrm{d}#2}}
\newcommand{\D}{\partial D}
\renewcommand{\S}{\mathcal{S}}
\newcommand{\outside}{\mathbb{R}^2\setminus \overline{D}}
\newcommand{\closure}[2][3]{%
	{}\mkern#1mu\overline{\mkern-#1mu#2}}
\renewcommand*{\Re}{\operatorname{Re}}
\renewcommand*{\Im}{\operatorname{Im}}

\newcommand{\de}{\: \mathrm{d}}

\renewcommand{\a}{\boldsymbol{\alpha}}

\newcommand{\N}{\mathbf{N}}

\newcommand{\F}{\mathcal{F}}
\newcommand{\R}{\tilde{R}_k}
\newcommand{\tpsi}{\tilde{\psi}_k}

\numberwithin{equation}{section}

\graphicspath{ {./figures/} }

\usepackage{cite}

\definecolor{ETHc}{RGB}{18,105,176}
\definecolor{blizzardblue}{rgb}{0.67, 0.9, 0.93}

\title{Mimicking the active cochlea with a fluid-coupled array of subwavelength Hopf resonators}
\author{Habib Ammari\thanks{Department of Mathematics, ETH Z\"urich, R\"amistrasse 101, CH-8092 Z\"urich, Switzerland (habib.ammari@math.ethz.ch, bryn.davies@sam.math.ethz.ch).} \qquad Bryn Davies\footnotemark[1]}
\date{}

\begin{document}
	
	\maketitle
	\noindent\rule{\textwidth}{0.4pt}
	\begin{abstract}
		We present a design for an acoustic metamaterial that mimics the behaviour of the active cochlea. This material is composed of a size-graded array of cylindrical subwavelength resonators, has similar dimensions to the cochlea and is able to reproduce the frequency separation of audible frequencies. Non-linear amplification is introduced to the model in order to replicate the behaviour of the cochlear amplifier. This formulation takes the form of a fluid-coupled array of Hopf resonators.	We seek solutions based on a modal decomposition, so as to retain the physically-derived coupling between resonators.
		
	\end{abstract}
	
	
	
	\noindent\rule{\textwidth}{0.4pt}
	



\section{Introduction}
\subsection{The cochlea}

The cochlea is the part of the ear responsible for transforming physical vibrations into neural signals. It does so with remarkable efficacy, capturing a large range of volumes with impressive frequency resolution. At its simplest, the cochlea is a long tube that is filled with fluid and divided in two by the cochlear partition \cite{dallos1996cochlea}. This partition is set in motion by acoustic waves in the cochlear fluid and contains receptor cells which respond to this motion \cite{hudspeth1983hair, reichenbach2014physics, martin2001compressive}. These cells, known as \emph{hair cells}, can be split into two categories, referred to as \emph{inner} and \emph{outer} due to their positioning. The former are responsible for converting motion of the partition into electrical signals, while the latter contain a motor protein and are involved in the cochlea's non-linear amplification mechanism.



The cochlea's mechanism for resolving frequency is based on the fact that several of its physical parameters (such as the stiffness and width of the basilar membrane, upon which the cochlear partition is mounted) are graded along its length. This property means that different tones will lead to maximal excitation at different points \cite{von1960experiments}. High frequency sounds are detected at the base of the cochlea while sounds of lower frequencies are detected closer to the apex. This frequency separation mechanism is remarkably effective, meaning that even non-musicians can detect differences of less than 0.5\% \cite{hudspeth2008making,dallos1992active}.

Humans are also able to comfortably hear sounds with amplitudes that range over an impressive six orders of magnitude. This is thanks to the cochlea's non-linear amplification response, known as the \emph{cochlear amplifier}. This is such that a compressive non-linearity occurs, whereby quieter sounds are amplified much more greatly than louder ones. Understanding the details of the mechanism responsible for the cochlear amplifier represents the most significant open problem in the study of hearing. Recently, the most prevalent models have been based on a critically-poised Hopf resonator, the details of which are discussed in \Cref{sec:hopf_intro}.


\subsection{An acoustic metamaterial}

Metamaterials are broadly defined as synthetic structures which have remarkable properties that do not generally exist in conventional materials. Typically composed of a repeating small structure, acoustic metamaterials have been successfully used for a variety of wave-control applications \cite{cummer2016controlling, craster2012acoustic, ma2016acoustic}. In this work we will study a gradient acoustic metamaterial, where the material parameters (in particular, the size of the small, repeating unit) vary with position. It is known that metamaterials of this kind can perform spatial frequency separation. This phenomenon is often known as \emph{rainbow trapping}, particularly in the context of optics where it was first proposed by Tsakmakidis et al. \cite{tsakmakidis2007trapped}. It has since been observed a in variety of different settings including acoustic metamaterials \cite{zhu2013acoustic, davies2019fully, jimenez2017rainbow, chen2014enhanced}, amongst others \cite{bennetts2018graded, skelton2018multi, jang2011plasmonic, colquitt2017seismic}.

%
In this work, we study an acoustic metamaterial composed of an array of high-contrast subwavelength resonators, graded in size. These are compressible objects with material parameters (density and bulk modulus) that differ greatly from the background medium and, as a result, have a resonant frequency that corresponds to a wavelength much greater than the object's physical size. The classical example is an air bubble in water (known as a \emph{Minnaert bubble} \cite{minnaert1933musical}) but the theory holds for any pair of materials with a sufficiently great parameter contrast. An asymptotic analysis of a graded array of such resonators was recently performed by Ammari and Davies \cite{davies2019fully}, where it was shown that with a suitably-chosen gradient the structure reproduced the spatial frequency separation exhibited by the cochlea \cite{von1960experiments}.

Rupin et al. \cite{rupin2019mimicking} have recently produced promising experimental results, demonstrating that a graded array of cylindrical resonators can mimic both the cochlea's frequency separation and non-linear amplification properties. Their set-up was scaled up for practicality and employed an array of cylindrical quarter-wavelength resonators with a microphone-and-speaker amplification system. In this work we study a structure with similar dimensions to the cochlea, consisting of gas-filled cylinders surrounded by fluid. We take a cross section of this structure so that we study the problem of an acoustic wave scattered in a plane by a size-graded array of circles.

\begin{figure}
	\begin{center}
		\begin{tikzpicture}[scale=0.5]
		\draw[fill=gray!30!white] (0,0.7) -- (0,1) -- plot [smooth] coordinates {(0,1) (9,0.9) (10.5,0) (9,-0.9) (0,-1)} -- (0,0.2) -- plot [smooth] coordinates {(0,0.2) (0.1,0.45) (0,0.7)};
		\draw plot [smooth] coordinates {(0,0.2) (-0.1,0.45) (0,0.7)};
		\draw[line width=1.5pt] (0,0) -- (9,0);
		
		\node at (-0.8,0.45) {\tiny OW};
		\node at (4.5,0.3) {\tiny CP};
		\node at (9.9,1.4) {\tiny Apex};
		\node at (0.3,1.4) {\tiny Base};
		
%
		
		
		\draw[white, fill = white, semithick] 
		(5,1.3) rectangle (9,5) (7.05,3) circle (1.5);
		\draw[gray] (7.05,0) circle (0.1);
		\draw[gray] (6.96,-0.05) -- (5.78,2.2);
		\draw[gray] (7.14,-0.05) -- (8.32,2.2);
		\draw[gray] (7.05,3) circle (1.5);
		
		\begin{scope}[xshift=4cm,yshift=-0.2cm]
		\draw[fill=gray,opacity=0.6] (1.4,2.7) rectangle (4.7,3.75);
		\draw (1.4,3.2) -- (4.7,3.2);
		
		\draw (3.1,2.95) circle (0.15);
		\draw (3.23,2.87) -- (3.51,3.42);
		\draw (2.96,3) -- (3.24,3.55);
		\draw[domain=-30:160] plot ({3.38+0.15*cos(\x)}, {3.5+0.15*sin(\x)});
		
		\begin{scope}[xshift=0.6cm]
		\draw (3.1,2.95) circle (0.15);
		\draw (3.23,2.87) -- (3.51,3.42);
		\draw (2.96,3) -- (3.24,3.55);
		\draw[domain=-30:160] plot ({3.38+0.15*cos(\x)}, {3.5+0.15*sin(\x)});
		\end{scope}
		
		\begin{scope}[xshift=1.2cm]
		\draw (3.1,2.95) circle (0.15);
		\draw (3.23,2.87) -- (3.51,3.42);
		\draw (2.96,3) -- (3.24,3.55);
		\draw[domain=-30:160] plot ({3.38+0.15*cos(\x)}, {3.5+0.15*sin(\x)});
		\end{scope}
		
		\begin{scope}[xshift=-0.6cm]
		\draw (3.1,2.95) circle (0.15);
		\draw (3.23,2.87) -- (3.51,3.42);
		\draw (2.96,3) -- (3.24,3.55);
		\draw[domain=-30:160] plot ({3.38+0.15*cos(\x)}, {3.5+0.15*sin(\x)});
		\end{scope}
		
		\begin{scope}[xshift=-1.2cm]
		\draw (3.1,2.95) circle (0.15);
		\draw (3.23,2.87) -- (3.51,3.42);
		\draw (2.96,3) -- (3.24,3.55);
		\draw[domain=-30:160] plot ({3.38+0.15*cos(\x)}, {3.5+0.15*sin(\x)});
		\end{scope}
		
		\begin{scope}[xshift=-1.8cm]
		\draw (3.1,2.95) circle (0.15);
		\draw (3.23,2.87) -- (3.51,3.42);
		\draw (2.96,3) -- (3.24,3.55);
		\draw[domain=-30:160] plot ({3.38+0.15*cos(\x)}, {3.5+0.15*sin(\x)});
		\end{scope}
		
		\node at (3.25,2.4) {\tiny BM};
		\end{scope}
		
		\draw[white, fill = white, semithick] 
		(5,1.3) rectangle (9,5) (7.05,3) circle (1.5);
		\draw[gray] (7.05,0) circle (0.1);
		\draw[gray] (6.96,-0.05) -- (5.78,2.2);
		\draw[gray] (7.14,-0.05) -- (8.32,2.2);
		\draw[gray] (7.05,3) circle (1.5);
		
		\begin{scope}[xshift=14cm,scale=1.3]
		\node at (7,0.8) {\tiny Resonator array};
		
		\foreach \x in {0,1,...,9} {
			\draw[fill=gray!30!white] (0.2*30*1.1^\x-0.2*30+0.1,0) circle (0.2*1.1^\x);}
		
		\begin{scope}[xshift=0cm,yshift=1cm]
		\draw [->, thick, opacity=0.7] (-0.1,0) -- (0.5,0);
		\draw [->, thick, opacity=0.7] (0,-0.1) -- (0,0.5);
		\node[opacity=0.7] at (0.8,-0.05) {\small $x_1$};
		\node[opacity=0.7] at (0,0.7) {\small $x_2$};
		\end{scope}
		\end{scope}
		\end{tikzpicture}
	\end{center}
	\caption{\textit{Left:}~A cross-section of a simplified, straightened-out model of the cochlea. The cochlea is partitioned along its length by the cochlear partition (CP). Signals enter the upper channel through the oval window (OW) and cause the partition to vibrate. \textit{Inset:}~The basilar membrane (BM), upon which the CP is mounted, is spanned radially by collagen fibres. \textit{Right:}~We mimic the cochlea with an array of circular resonators which are size-graded to replicate the properties of the BM.} \label{fig:cochlea_diagram}
\end{figure}
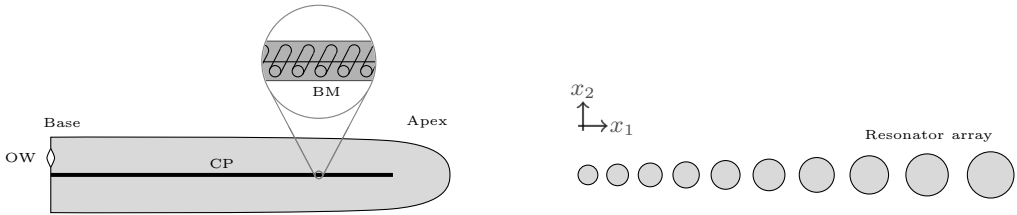

As an aside, a parallel source of inspiration for the gradient metameterial considered in this work is the assortment of cochlear models based on an arrays of harmonic oscillators \cite{babbs2011quantitative, wada1998measurement, lerud2019canonical, kern2003essential, wilkinson1921, wilson1992cochlear}. The compliance of the cochlear partition is largely accounted for by the basilar membrane which is composed of soft, elastic tissue reinforced by strong collagen fibres across the cochlea's width (\Cref{fig:cochlea_diagram}) \cite{liu2008orthotropic, naidu2007basilar, reichenbach2014physics}. As a result, if one considers the membrane in sections then the mechanical coupling between them is relatively small and each can be modelled as an individual harmonic oscillator. This motivates discrete cochlear models based on graded arrays of, for example, masses on springs \cite{babbs2011quantitative} or elastic beams \cite{wada1998measurement}. In fact, many of the earliest theories of hearing were based on a version of this principle \cite{bell2004resonance}, such as Helmholtz' famous ``piano strings'' analogy \cite{helmholtz1875sensations}. Since Minnaert-type resonators behave as harmonic oscillators \cite{devaud2008minnaert,minnaert1933musical}, the structure studied here can be viewed as a manifestation of these models. 

The choice of resonators made here is significant since it allows for the design of a structure that has similar dimensions to a cochlea and responds to audible frequencies. Throughout this paper we will run simulations on an array of 22 resonators, arranged in a linear configuration that measures $32\si{\milli\metre}$ in length and $0.28\si{\milli\metre}$ at its widest. We use the material parameters of air and water for the resonators and surrounding fluid, respectively, in order to demonstrate that the material properties (in particular, the contrast) required to give the desired subwavelength resonant response are exhibited by everyday substances. Bubbly structures of this kind have previously been constructed \emph{e.g.} by injecting air bubbles into silicone-based polymers \cite{leroy2009design, leroy2009transmission}.

Modelling the interactions between the array of resonators, which occur on subwavelength scales, is a challenging problem. We will use a modal-decomposition approach to study the wave-propagation problem \cite{wagg2016nonlinear, fletcher1992acoustic}. The structure's spatial eigenmodes will be found using layer-potential formulations \cite{ammari2009layer} before these profiles are used as a basis to study the behaviour of the system when appropriate non-linear amplification is introduced to the resonators.

\subsection{Hopf resonators in cochlear mechanics} \label{sec:hopf_intro}

While there is still much to be learnt about the cochlear amplifier, researchers are guided by a handful of valuable observations. It is understood, firstly, that this mechanism acts as a negative-damping mechanism \cite{gold1948hearing,neely1986model,joyce2017study}. Secondly, it has been observed that non-linear amplification occurs with a one-third power law, suggesting the presence of a cubic non-linearity \cite{hudspeth2008making}. 


Hopf resonators have become popular objects to study in the field of non-linear cochlear mechanics thanks to their remarkable ability to account for the key properties that typify cochlear behaviour \cite{eguiluz2000essential, hudspeth2008making, hudspeth2010critique, camalet2000auditory, joyce2017study, duke2008critical, kern2003essential, lerud2019canonical, magnasco2003wave}. The normal form of a single Hopf resonator $z=z(t):\mathbb{R}\to\mathbb{C}$ in the complex plane is given by the forced differential equation

\begin{equation} \label{eq:hopf}
\dd{z}{t} = (\mu+i\omega_0)z-|z|^2z + F,
\end{equation}
where $F=F(t)$ is the forcing term and $\omega_0$ and $\mu$ are real parameters. This system is a resonator in the sense that the absolute value of the response $z$ is greatest when the forcing $F$ occurs with frequency $\omega_0$. In cochlear models, $z$ is some variable which characterises the system's state.

The parameter $\mu$ is the bifurcation parameter. For $\mu<0$ the unforced system ($F=0$) has a stable equilibrium at $z=0$ whereas when $\mu>0$ this equilibrium is unstable and there exists a stable limit cycle given by $z(t)=\sqrt{\mu}e^{i\omega_0 t}$. This birth of a limit cycle is typical of a (supercritical) Hopf bifurcation, which is formally characterised by a conjugate-pair of linearised eigenvalues crossing the imaginary axis. Writing the unforced system \eqref{eq:hopf} in terms of its real and imaginary parts and linearising about the fixed point at zero gives a system whose Jacobian matrix has eigenvalues $\lambda=\mu\pm i\omega_0$. These eigenvalues clearly cross the imaginary axis when $\mu$ passes zero. For further details see \emph{e.g.} \cite{strogatz2018nonlinear, jordan1999nonlinear, wagg2016nonlinear}. 

The greatly enhanced response for frequencies close to $\omega_0$ is able to account for the cochlea's frequency selectivity. The cubic non-linearity in \eqref{eq:hopf} is able to reproduce the one-third power law of the cochlea: when $\mu$ is small, and the system is close to bifurcation, we have that $|z|\approx |F|^{1/3}$ for frequencies close to resonance.

One of the earliest pieces of evidence supporting the active nature of the cochlea was the observation that the ear emits sounds (known at otoacoustic emissions) as part of its response \cite{kemp2002otoacoustic, zurek1985acoustic}. The existence of stable limit cycles, for certain parameter values, predicts this behaviour \cite{camalet2000auditory, duke2008critical, kanders2019}.




A further symptom of the non-linearity that exists in the cochlea is the behaviour that is observed under the influence of a signal composed of two distinct tones. Firstly, when the ear is excited by such a stimulus two-tone suppression occurs. That is, the frequency spectrum of the response contains the expected two amplitude peaks, however, these are smaller than each would be in the absence of the other tone \cite{ruggero1992two}. Further, in this situation the ear also detects additional tones, variously known as combination tones, distortion products or Tartini's tones \cite{helmholtz1875sensations, robles1997two, julicher2001physical}. Close to bifurcation, the non-linearity in \eqref{eq:hopf} gives products that can account for these phenomena \cite{julicher2001physical, duke2008critical}.

In this work, we will introduce a Hopf-type non-linearity directly to the wave-propagation problem by supposing that the resonators are equipped with an appropriate forcing mechanism (as was realised by \cite{rupin2019mimicking}). Based on an eigenmode decomposition, we will explore the Hopf-type behaviour of this system and show that the crucial cochlea-like properties of \eqref{eq:hopf} are retained by the coupled subwavelength structure.

\section{Scattering by coupled subwavelength resonators} \label{sec:passive_coupling}

It is known that the frequency-location (tonotopic) map in the cochlea is exponential \cite{von1960experiments} and it was shown in \cite{davies2019fully} that if an array of subwavelength resonators has a similar size-grading, we reproduce this phenomenon. Thus, we will consider a domain $D$ in $\mathbb{R}^2$ which is the disjoint union of $N\in\mathbb{N}$ circular subdomains $\{D_1,\ldots,D_N\}$ with each successive radius defined as $R_{n+1}=sR_n$, for some $s>1$.

We denote by $\rho$ and $\kappa$ the density and bulk modulus of the interior of the resonators, respectively, and use $\rho_0$ and $\kappa_0$ for the corresponding parameters for the background fluid (which occupies $\mathbb{R}^2\setminus \overline{D}$). We may then denote the acoustic wave speeds in $\mathbb{R}^2\setminus \overline{D}$ and $D$, respectively, by
\begin{equation}
v_0=\sqrt{\tfrac{\kappa_0}{\rho_0}}, \quad v=\sqrt{\tfrac{\kappa}{\rho}}.
\end{equation}
We also introduce the dimensionless contrast parameter
\begin{equation} \label{defn:contrasts}
\delta := \frac{\rho}{\rho_0}.
\end{equation}
Since we want our structure to be of similar dimensions to the cochlea and exhibit a resonant response to audible frequencies, we need subwavelength resonant modes to exist. This is known to occur in the case that $\delta \ll 1$ \cite{minnaert1933musical, ammari2018minnaert, davies2019fully}. For the simulations performed here we use material parameters corresponding to air-filled resonators surrounded by water, giving that $\delta\approx 10^{-3}$.

We will use a modal decomposition to analyse the wave propagation within the structure \cite{wagg2016nonlinear,fletcher1992acoustic}. That is, we wish to express the acoustic pressure $p=p(x,t)$ at position $x$ and time $t$ in the form
\begin{equation} \label{eq:modal_decomp}
p(x,t) = \Re\left(\sum_{n}\alpha_n(t)u_n(x)\right).
\end{equation}
Separating variables in the unforced, linear wave equation yields a spatial Helmholtz problem given by

\begin{equation} \label{eq:helmholtz_equation}
\begin{cases}
\left( \Delta + \tfrac{\omega^2}{v_0^2} \right) u(x,\omega) = 0, & \text{for }(x,\omega)\in \outside\times\mathbb{C}, \\
\left( \Delta + \tfrac{\omega^2}{v^2} \right) u(x,\omega) = 0, & \text{for }(x,\omega)\in D\times\mathbb{C}, \\
u_+ - u_- = 0, & \text{for }(x,\omega)\in \D\times\mathbb{C},\\
\delta \ddp{u}{\nu}\big|_+ -  \ddp{u}{\nu}\big|_- = 0, & \text{for }(x,\omega)\in \D\times\mathbb{C},
\end{cases}
\end{equation}
where $\omega\in\mathbb{C}$ is the separation constant, $\ddp{}{\nu}$ denotes the outward normal derivative and the subscripts + and - are used to denote evaluation from outside and inside $\D$, respectively.
We must also insist that $u(\cdot,\omega)$ satisfies the Sommerfeld radiation condition
\begin{equation} \label{eq:src}
\lim\limits_{|x|\to\infty} |x|^{1/2} \left(\ddp{}{|x|}-i\frac{\omega}{v_0}\right)u(x,\omega)=0.
\end{equation}
This condition is required to ensure that the solution represents outgoing waves (rather than incoming from infinity) and gives the well-posedness of \eqref{eq:helmholtz_equation}. Since we study the problem in unbounded space, energy is lost to the far-field.

In light of the fact that \eqref{eq:helmholtz_equation} contains no forcing or amplification, we define a \emph{resonant frequency} and associated {\emph{eigenmode}} (or \emph{resonant mode}) to be solutions $(\omega,u(\cdot,\omega))\in\mathbb{C}\times H_{loc}^1(\mathbb{R}^2)$ of \eqref{eq:helmholtz_equation}. Here, $H_{loc}^1(\mathbb{R}^2)$ is the space of functions that, on every compact subset of $\mathbb{R}^2$, are square integrable and have a weak first derivative that is also square integrable. We are interested in solutions where $\omega$ is small and the resonators are much smaller than the wavelength of the associated radiation, such solutions are referred to as \emph{subwavelength} modes.

\subsection{Layer-potential approach}

We solve \eqref{eq:helmholtz_equation} using a layer-potential approach \cite{ammari2009layer}. This is based on representing the solution, in terms of some surface potentials $\phi,\psi\in L^2(\D)$, as
\begin{equation} \label{eq:layer_potential_representation}
	u(x,\omega) = \begin{cases}
	\S_{D}^{\omega/v_0}[\psi](x), & (x,\omega)\in\outside\times\mathbb{C},\\
	\S_{D}^{\omega/v}[\phi](x), & (x,\omega)\in D\times\mathbb{C},
	\end{cases}
\end{equation}
where $\S_D^\omega$ is the \emph{Helmholtz single layer potential} associated with the domain $D$. This integral operator is defined as
\begin{equation}
	\S_{D}^\omega[\varphi](x) := \int_{\D} \Gamma^\omega(x-y) \varphi(y)\, d\sigma(y),\quad x\in\D,\, \varphi\in L^2(\D),\, \omega\in\mathbb{C},
\end{equation}
where $\Gamma^\omega$ is the outgoing (\emph{i.e.} satisfying the Sommerfeld radiation condition) fundamental solution to the Helmholtz operator $\Delta+\omega^2$ in $\mathbb{R}^2$ \cite{ammari2018mathematical}. The value of this approach is that solving \eqref{eq:helmholtz_equation} is reduced to finding $\phi,\psi$ such that the two transmission conditions on $\D$ hold \cite{ammari2009layer,ammari2018minnaert}.

An asymptotic analysis of the resonant frequencies and eigenmodes, based on their layer-potential representations \eqref{eq:layer_potential_representation}, was performed in \cite{davies2019fully}. It was shown that a system of $N$ coupled resonators has $N$ subwavelength resonant modes $u_1(x),\dots,u_N(x)$ and corresponding resonant frequencies $\omega_1,\dots,\omega_N$ with positive real part.

Using the layer-potential representation \eqref{eq:layer_potential_representation}, we are able to find the subwavelength resonant modes numerically by expanding the functions $\phi$ and $\psi$ in terms of Fourier bases on the boundary of each resonator. This is particularly convenient in the case of circular resonators since the evaluation of $\S_D^\omega$ on Fourier modes has closed-form expressions, the details of which are given in \Cref{app:multipole}.

\begin{figure}[t]
	\begin{center}
		\includegraphics[width=\textwidth,trim={2cm 2cm 0 1.5cm}]{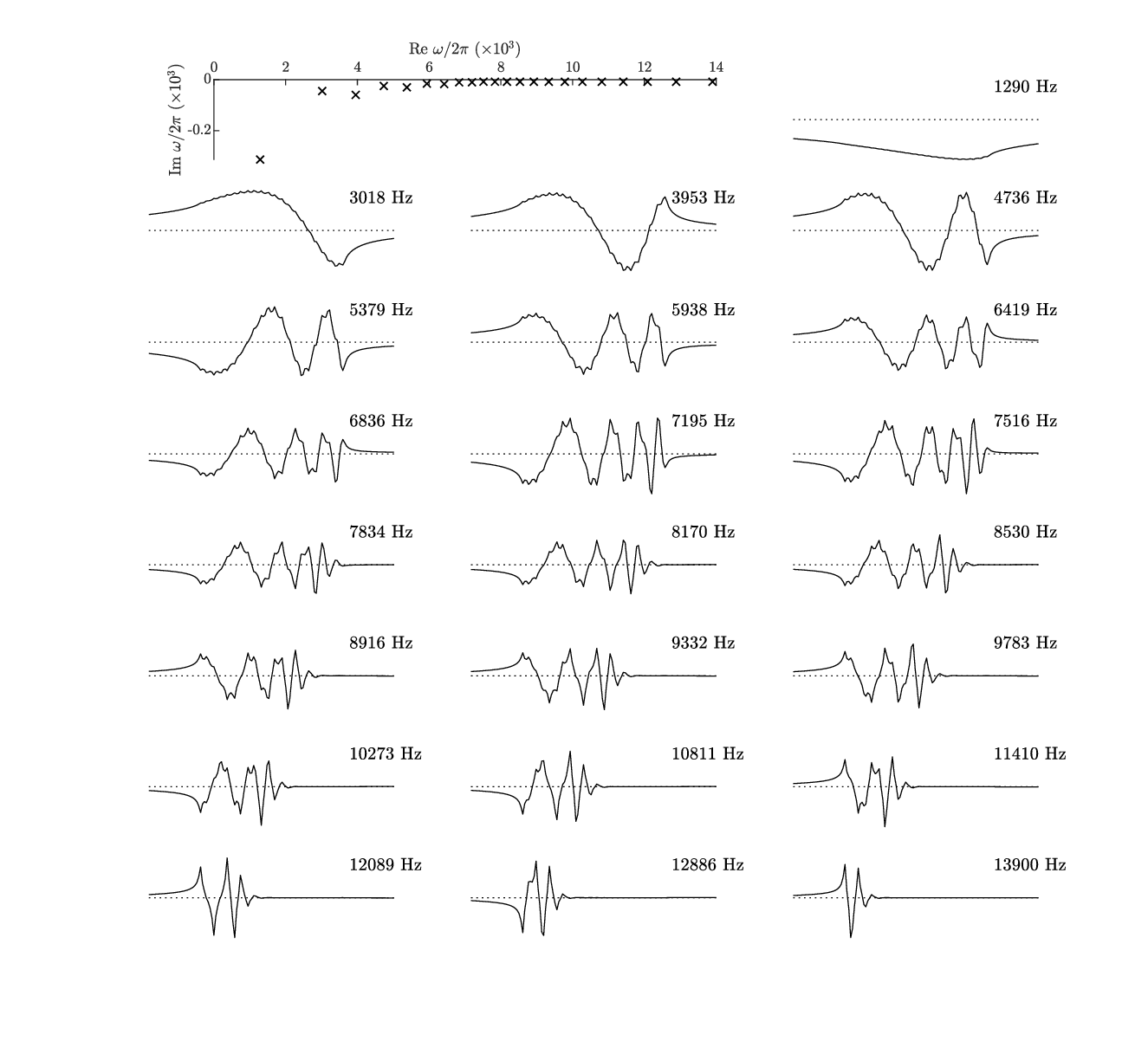}
		\caption{The eigenmodes and associated resonant frequencies for the system of 22 resonators. The eigenmodes are plotted along the line through the resonators' centres and the resonant frequencies are shown in the complex plane.} \label{fig:modes}
	\end{center}
\end{figure}

The $N$ subwavelength eigenmodes take the form of increasingly oscillating coupled patterns (\Cref{fig:modes}). The profiles share several similarities with the response of the basilar membrane. Each mode has a position of maximal amplitude (sometimes known as the resonant place) that depends on its resonant frequency \cite{von1960experiments, davies2019fully}. In the region of the resonant place, the wavelength decreases as the amplitude peaks, before the solution dies away quickly. This was similarly observed by \cite{rupin2019mimicking} in simulations using point scatterers. The use of resonators with non-zero radii, and the fact that each mode is approximately constant on each resonator \cite{davies2019fully}, mean the profiles look less smooth here.

%
%

The design studied here (with $R_1=0.1\si{\milli\metre}$ and $s=1.05$) is chosen so that resonant frequencies (\Cref{fig:modes}) have real parts which fall within the range of audible frequencies (often quoted as 20~Hz~-~20~kHz). The negative imaginary parts denote the positive rate of attenuation associated with each mode, due to energy being lost to the far field (note that the assumed time harmonicity in \eqref{eq:src} is $e^{-i\omega t}$). We will see in \Cref{sec:nonlinear} that the amplification required to illicit a Hopf bifurcation depends on this quantity.

The structure will also have higher order resonant modes at frequencies which correspond to wavelengths similar to the size of the resonators, or bigger. Given the physical dimensions and wavelengths of the problem we are interested in, we focus our attention on the subwavelength modes as these will dominate the behaviour of the system.


\section{Non-linear system} \label{sec:nonlinear}

We now wish to introduce appropriate non-linear amplification to the model. As discussed in \Cref{sec:hopf_intro}, the canonical form of a Hopf resonator is able to account for the important properties of the cochlear amplifier. This suggests adding amplification based on a non-linearity of the form
\begin{equation} \label{eq:nonlinearity}
\mathcal{N}[\varphi]:=\mu \varphi-\beta|\varphi|^2\varphi.
\end{equation}
These two terms, respectively, account for the negative damping and cubic non-linearity that we said our amplification should include.

An important consideration, when choosing to introduce amplification, is the stability of the system. For example, Rupin et al. \cite{rupin2019mimicking} used a formulation whereby amplification closely resembling $\mathcal{N}[p]$ was added. In order for this formulation to be stable, it was necessary to design a set-up that switched off the amplification if the pressure exceeded a threshold value. Conversely, there exist a number of formulations which are stable without this thresholding. Examples include variants of $\mathcal{N}[\partial_t p]$ used in the artificial cochlear devices of Joyce and Tarazaga \cite{joyce2017study, joyce2014mimicking, joyce2015developing} and the $\mu \partial_t p-\beta|p|^2\partial_t p$ term considered by Duke and J{\"u}licher \cite{duke2008critical}. 

In this work we will study the system produced by introducing amplification of the form $\mathcal{N}[\partial_t p]$ to the resonators. As we shall see in \Cref{sec:stability}, this system is stable without the need to impose a pressure threshold. Furthermore, there is evidence which suggests that hair cell stimulation (by stereocilia displacement) is dependent on membrane velocity \cite{lu2009stereocilia}, suggesting that any amplification should be a function of $\partial_t p$.



Since the $N$ subwavelength modes are expected to dominate, we truncate the expansion \eqref{eq:modal_decomp}, seeking a solution of the form
\begin{equation} \label{eq:soln_decomposition}
p(x,t) = \Re\left(\sum_{n=1}^{N} \alpha_n(t)u_n(x)\right),
\end{equation}
for some complex-valued functions of time $\alpha_1(t),\dots,\alpha_N(t)$. We will also project the forcing term onto the space spanned by the subwavelength eigenmodes. If $Q\subset\mathbb{R}^2$ is a compact set on which the forcing is applied (and $\mathcal{X}_Q$ is the characteristic function of $Q$) then we decompose the forcing as
\begin{equation} \label{eq:forcing}
f(t)\mathcal{X}_Q(x)\simeq f(t)\sum_{n=1}^{N} \F_nu_n(x),
\end{equation}
where $\F_n:=(\mathcal{X}_Q,u_n)_{2,Q}$.

In light of the transmission properties (across $\D$) that the eigenmodes inherit from \eqref{eq:helmholtz_equation}, we reach the problem
\begin{equation} \label{eq:decomposed_equation}
\sum_{n=1}^{N} \big(\alpha_n''(t)+\omega_n^2\alpha_n(t)\big)u_n(x) =f(t)\sum_{n=1}^{N} \F_nu_n(x) + \mathcal{N}\left[\sum_{n=1}^{N} \alpha_n'(t)u_n(x)\right] \mathcal{X}_D(x).
\end{equation}


\subsection{Modal system} \label{sec:normal_form}

Our approach to studying \eqref{eq:decomposed_equation} will be to take the $L^2(D)$ product with $u_m$ for $m=1,\dots,N$ to reach a coupled system of $N$ ordinary differential equations. Define the matrix $\gamma\in\mathbb{C}^{N\times N}$ as $\gamma_{ij}:=(u_i,u_j)_{L^2(D)}$. 
An important property is that, thanks to the linear independence of the eigenmodes, $\gamma$ is invertible \cite{davies2019fully}.
If we write $\alpha_1,\dots,\alpha_N$ in the column vector $\a$ then the modal system is described by
\begin{equation} \label{eq:coupled_syst}
\a''+\mathbf{\Lambda}\a-{\mu}\a' +\beta\N(\a')  = f(t)\mathbf{F},
\end{equation}
where $\mathbf{\Lambda}\in\mathbb{C}^{N\times N}$ is diagonal with entries $\omega_n^2$, $\mathbf{F}\in\mathbb{C}^{N}$ is the vector of forcing constants $\F_n$ and $\N:\mathbb{C}^{N}\to\mathbb{C}^{N}$ is the non-linear function defined as
\begin{equation} \label{eq:coupled_nonlinearity}
\N(\mathbf{z}) := \gamma^{-1}\left[\bigintsss_D\left|\sum_{n=1}^{N} z_nu_n(x)\right|^2 \left(\sum_{n=1}^{N} z_nu_n(x)\right)\closure{u_j(x)}\de x \right]_{j=1,...,N}.
\end{equation}

\begin{figure}[t]
	\begin{center}
		\includegraphics[width=0.9\textwidth]{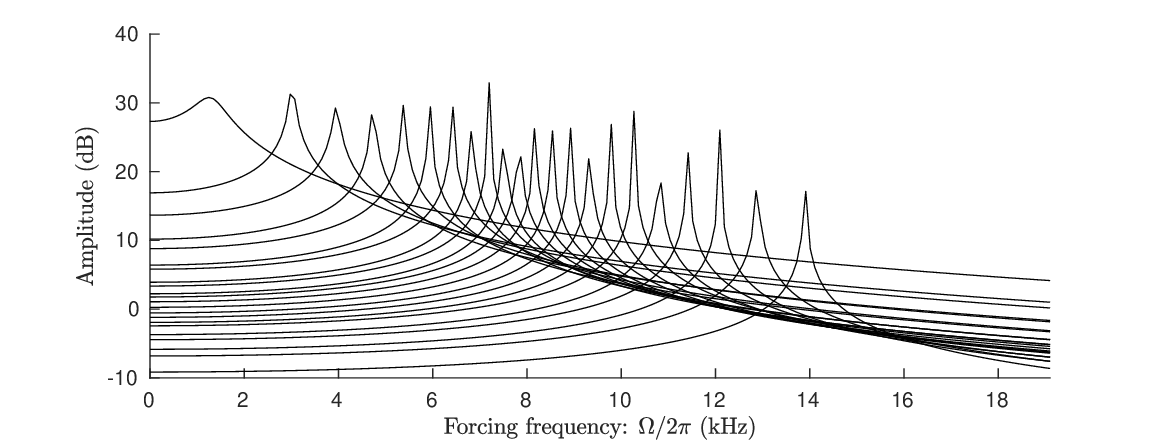}
		\caption{The response of each amplitude $|X_1|,\dots,|X_N|$ (with peaks from left to right) as a function of the incident frequency $\Omega$. Each mode experiences a peak of excitation in the vicinity of its resonant frequency. The incident sound is at 100~dB~SPL.} \label{fig:full_decomp}
	\end{center}
\end{figure}

In much of what follows we will be interested in the case when $f(t)$ is harmonic with frequency $\Omega$. In this case, we can approximate the solution to \eqref{eq:coupled_syst} using a harmonic balance approach \cite{wagg2016nonlinear,stoker1950nonlinear,fletcher1992acoustic,nayfeh1980nonlinear}. That is, if $f(t)=Fe^{-i\Omega t}$, for $F,\Omega\in\mathbb{R}$, then we may approximate the steady-state solutions to \eqref{eq:coupled_syst} as $\alpha_k(t)=X_k e^{-i\Omega t+i\psi_k}$ for amplitudes $X_k\in\mathbb{R}$ and phase delays $\psi_k\in\mathbb{R}$. Making this substitution leads to a system of coupled cubic equations that can be solved numerically. 

In \Cref{fig:full_decomp} we see that, as is to be expected, that as the forcing frequency is varied each mode is excited much more greatly in the vicinity of the associated resonant frequency (in spite of the coupling within the non-linearity \eqref{eq:coupled_nonlinearity}). This motivates an approximate system whereby, if the system is forced at a frequency close to one of the resonant frequencies, we assume that only that mode is excited.




\subsection{Single-mode approximation} \label{sec:single_mode}

When the forcing frequency $\Omega$ is close to one of the resonant frequencies $\omega_k$ we approximate the solution to \eqref{eq:decomposed_equation} by assuming that only the corresponding mode $u_k$ is excited. In such a regime, we take the $L^2(D)$ product of \eqref{eq:decomposed_equation} with $u_k$ to yield the equation
\begin{equation} \label{eq:ode}
\alpha_k''+ \omega_k^2 \alpha_k- \mu\alpha_k'+ \hat{\beta}|\alpha_k'|^2\alpha_k'= f(t)\F_k, 
\end{equation}
where $\hat{\beta}:=\beta\|u_k\|_{4,D}^4/\|u_k\|_{2,D}^2$.

\subsubsection{Hopf bifurcation} \label{sec:hopf}

%
%
%

At this point, we pause to explore the Hopf-type behaviour that is exhibited by our model. In the case that $f=0$, we see that \eqref{eq:ode} has a periodic solution $\alpha_k(t) = R_k^ce^{-i\Omega_k^ct}$ provided that $\mu\geq \mu_k^c$, where

\begin{equation} \label{eq:limit_cycle}
\begin{gathered}
\Omega_k^c:=\sqrt{\Re(\omega_k)^2-\Im(\omega_k)^2},\\ \mu_k^c:=\frac{-2\Re(\omega_k)\Im(\omega_k)}{\Omega_k^c},\qquad
R_k^c=\sqrt{\frac{\mu-\mu_k^c}{\hat{\beta}}}\frac{1}{\Omega_k^c}.
\end{gathered}
\end{equation}
This birth of a limit cycle is typical of a Hopf bifurcation. A Hopf bifurcation is characterised by a conjugate pair of linearised eigenvalues crossing the imaginary axis \cite{strogatz2018nonlinear}. Decomposing $\alpha_k$ into its real and imaginary parts, we can write \eqref{eq:ode} as a four-dimensional system of first-order ordinary differential equations. Linearising this system around the fixed point at $\alpha_k=0$ gives the Jacobian matrix
\begin{equation} \label{eq:J}
J = \begin{bmatrix}
0 & 0 & 1 & 0 \\
0&0&0&1\\
-(\Omega_k^c)^2 & \mu_k^c\Omega_k^c & \mu & 0\\
-\mu_k^c\Omega_k^c &-(\Omega_k^c)^2 & 0&\mu
\end{bmatrix},
\end{equation}
which has eigenvalues given, using the notation of \eqref{eq:limit_cycle}, by
\begin{equation}
\lambda=\frac{1}{2}\left(\mu\pm\sqrt{-4(\Omega_k^c)^2\pm 4i\Omega_k^c\mu_k^c+\mu^2}\right).
\end{equation}
When $\mu=\mu_k^c$, the eigenvalues of $J$ are $\lambda=\pm i\Omega_k^c,\mu_k^c\pm i\Omega_k^c$. It can also be shown that $\frac{\de}{\de \mu}\Re(\lambda)>0$ meaning that the pair of eigenvalues cross the imaginary axis (from left to right) as $\mu$ passes the critical value. 

\begin{figure}
	\begin{center}
		\includegraphics[width=0.8\textwidth,trim={0.5cm 0 0 0}]{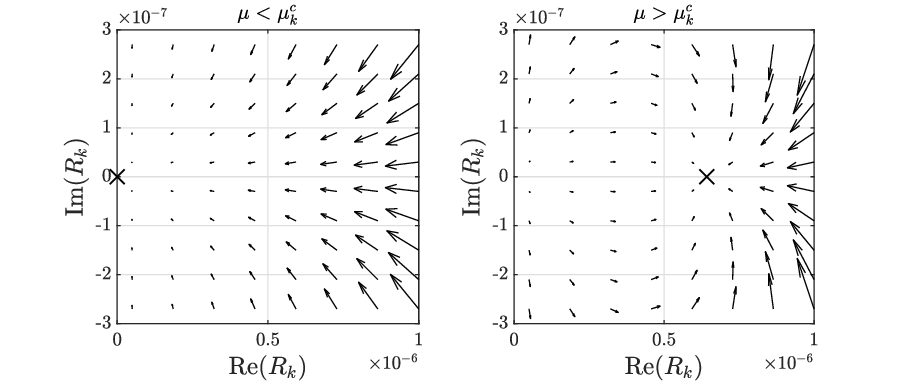}
		\caption{The birth of a limit cycle at Hopf bifurcation. For $\mu<\mu_k^c$ the origin is a stable equilibrium. When $\mu>\mu_k^c$ a stable limit cycle with radius $R_k^c$ is born. We depict $\mu=\mu_k^c\pm100$ and show the stable equilibria with crosses. Here, $k=11$ and $\hat{\beta}=10^5\si{\second\per\pascal\squared}$.} \label{fig:stability}
	\end{center}
\end{figure}

In order to visualise the local stability of this limit cycle and the fixed point at $\alpha_k=0$, we allow the radius $R_k$ to be a slowly varying function of $t$. That is, we use the ansatz $\alpha_k=R_k(t)e^{-i\Omega_k^ct}$ in \eqref{eq:ode} and disregard any terms containing either $R_k''$ or products of $R_k'$. This approach leads to the equation
\begin{equation} \label{eq:linear_stability}
(-2i\Omega_k^cR_k'-(\Omega_k^c)^2R_k+ \omega_k^2R_k)- \mu(R_k'-i\Omega_k^cR_k)+\hat{\beta} ((\Omega_k^c)^2R_k^2R_k'-i(\Omega_k^c)^3R_k^3)=0,
\end{equation}
which is linear in $R_k'$. The phase planes for $\mu>\mu_k^c$ and $\mu<\mu_k^c$ (\Cref{fig:stability}) demonstrate that as $\mu$ passes the critical value a stable limit cycle (described by \eqref{eq:limit_cycle}) is born out of the stable equilibrium at the origin, as is typical of a (supercritical) Hopf bifurcation.

\subsubsection{Stability of solutions} \label{sec:stability}

An important consideration, when choosing an appropriate non-linearity, is the stability of the system. We explore the stability of unforced solutions to \eqref{eq:ode} using a technique known as averaging \cite{wagg2016nonlinear, rand2012lecture, nayfeh1980nonlinear}. The analysis in \Cref{sec:hopf} showed that there is a locally stable limit cycle (when $\mu>\mu_k^c$) but it is valuable to understand what happens if $R_k(0)$ is further away from $R_k^c$. We begin with the ansatz
\begin{equation} \label{eq:average_ansatz}
\alpha_k(t)=R_k(t)e^{-i\Omega_k^c t+i\psi_k(t)}, \qquad \alpha_k'(t)=-i\Omega_k^c R_k(t)e^{-i\Omega_k^c t+i\psi_k(t)}.
\end{equation}
We assume for this analysis that $\hat{\beta}=:\epsilon>0$ and that $\mu=\mu_k^c+\epsilon$.

Differentiating the first expression of \eqref{eq:average_ansatz} and substituting into the second yields
\begin{equation} \label{eq:77}
\left(R_k'+i\psi_k'R_k\right)e^{-i\Omega t+i\psi_k}=0.
\end{equation}
Further, substituting \eqref{eq:average_ansatz} into \eqref{eq:ode} gives
\begin{equation} \label{eq:78}
\left(-i\Omega_k^c R_k' +\Omega_k^c\psi_k'R_k\right)e^{-i\Omega t+i\psi_k}=\epsilon i \left(-\Omega_k^cR_k+(\Omega_k^c)^3R_k^3\right) e^{-i\Omega t+i\psi_k}.
\end{equation}
We may take the real parts of \eqref{eq:77} and \eqref{eq:78} and solve for $R_k'$ and $\psi_k'$ to give
\begin{align}
R_k'&={\epsilon}\left(R_k-(\Omega_k^c)^2R_k^3\right)\sin^2(\Omega_k^c t-\psi_k), \label{eq:R_eqn}\\
\psi_k'&={\epsilon}\left(-R_k+(\Omega_k^c)^2R_k^3\right)\sin(\Omega_k^c t-\psi_k)\cos(\Omega_k^c t-\psi_k).
\end{align}

We now make a near-identity transformation in order to express $R_k$ and $\psi_k$ in terms of their average values over the interval $(t-\pi/\Omega_k^c,t+\pi/\Omega_k^c)$, which we denote by $\R$ and $\tpsi$. This transformation has the form
\begin{align} \label{eq:ni_transform}
R_k&=\R+\epsilon h_1(\R,\tpsi,t)+O(\epsilon^2), \\
\psi_k&=\tpsi+\epsilon h_2(\R,\tpsi,t)+O(\epsilon^2),
\end{align}
where $h_1$ and $h_2$ should be chosen in order to simplify the equations for $\R$ and $\tpsi$ as much as possible. This substitution leads to the equations
\begin{align}
\R'&={\epsilon}\left(-\ddp{h_1}{t}+\left(\R-(\Omega_k^c)^2\R^3\right)\sin^2(\Omega_k^c t-\tpsi)\right) + O(\epsilon^2), \label{eq:Rbar} \\
\tpsi'&=\epsilon\left(-\ddp{h_2}{t}+\left(-\R+(\Omega_k^c)^2\R^3\right)\sin(\Omega_k^c t-\tpsi)\cos(\Omega_k^c t-\tpsi)\right) + O(\epsilon^2).\label{eq:psibar}
\end{align}
Ideally, we would like to choose $h_1$ so that it cancels with the other $O(\epsilon)$ term in \eqref{eq:Rbar}. However, this antiderivative might grow in time meaning the expansion \eqref{eq:ni_transform} will not be valid for large $t$. Instead, we take $h_1$ as the antiderivative minus a linear term that grows with the average value \cite{rand2012lecture}, that is
\begin{equation} \label{eq:h1}
\begin{split}
h_1(\R,\tpsi,t)&=\int_{0}^{t} \left(\R-(\Omega_k^c)^2\R^3\right)\sin^2(\Omega_k^c t-\tpsi) \de t \\
&\qquad\qquad-\left[\frac{\Omega_k^c}{2\pi}\int_{0}^{2\pi/\Omega_k^c} \left(\R-(\Omega_k^c)^2\R^3\right)\sin^2(\Omega_k^c t-\tpsi) \de t\right]t.
%
\end{split}
\end{equation}

After substitution of \eqref{eq:h1} into \eqref{eq:Rbar}, we make an approximation in the spirit of the ``averaging'' methodology \cite{wagg2016nonlinear, rand2012lecture, nayfeh1980nonlinear}. We will assume that the integral in the second term of \eqref{eq:h1} can be well-approximated by taking the value of $\R$ and $\tpsi$ as constant over a cycle of oscillation, leaving a simple trigonometric integral.

We choose $h_2$ similarly and find that, up to an error of order $O(\epsilon^2)$,
\begin{equation}
\R'=\frac{1}{2}\epsilon\left(\R-(\Omega_k^c)^2\R^3\right),\qquad
\tpsi'=0.
\end{equation}
Solving by separation of variables gives that
\begin{equation}
R_k(t)=\frac{1}{\sqrt{R_k(0)^{-2}e^{-\epsilon t}+(\Omega_k^c)^2\left(1-e^{-\epsilon t}\right)}}+O(\epsilon).
\end{equation}

Crucially, for any $R_k(0)>0$ it holds that $R_k(t)\to R_k^c$ as $t\to\infty$, demonstrating that this limit cycle is asymptotically stable.

%
%

\subsubsection{Pure-tone response} \label{sec:amplification_selectivity}

\begin{figure}[t]
	\begin{center}
		\includegraphics[width=0.9\textwidth,trim={0 0 0 1cm}]{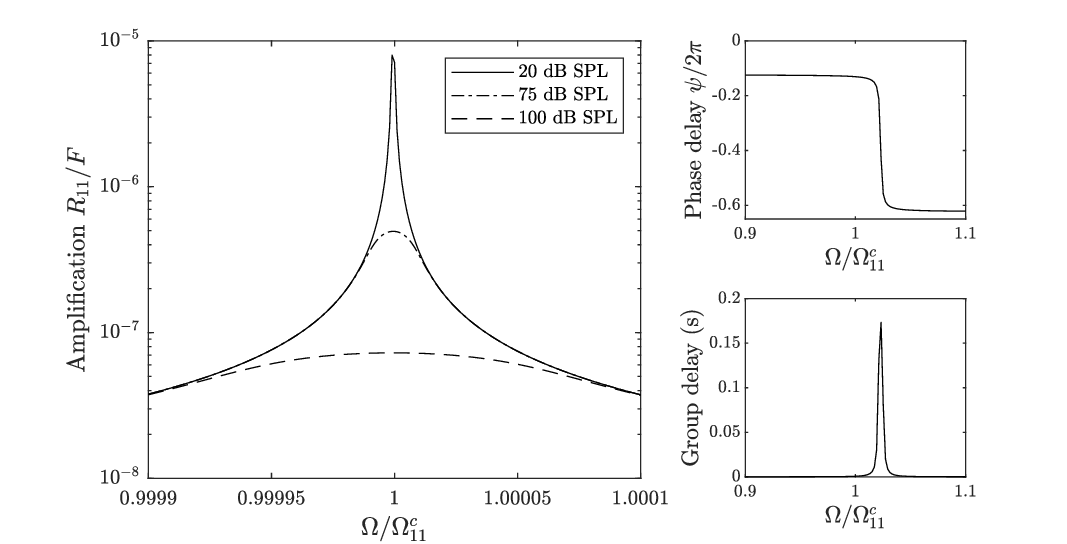}
		\caption{The non-linear response of the single-mode system. Amplification scales non-linearly with amplitude, and is greater for quieter sounds. Close to resonance, a phase delay of half a cycle is accumulated as well as a sharp increase in group delay. We take $\hat{\beta}=10^5\si{\second\per\pascal\squared}$ and $\mu=\mu_{11}^c$ so that the system is poised at bifurcation. The delay plots are shown for 20~dB~SPL.} \label{fig:ampli_delay}
	\end{center}
\end{figure}

Consider the case of an incoming signal that consists of a single pure tone at frequency $\Omega$, that is, $f(t)=Fe^{-i\Omega t}$ for $F,\Omega\in\mathbb{R}$, where $\Omega$ is close to $\omega_k$. Using the harmonic balance ansatz $\alpha_k(t)=R_k e^{-i\Omega t+i\psi_k}$ and finding the complex modulus of the resulting equation, we arrive at the amplitude-frequency response relation
\begin{equation} \label{eq:ampli_freq}
\left((\Omega_k^c)^2-\Omega^2\right)^2R_k^2 + \left(-\mu_k^c \Omega_k^c R_k+\mu\Omega R_k-\hat{\beta}\Omega^3R_k^3\right)^2 = F^2 |\mathcal{F}_k|^2.
\end{equation}
There is a sharply increased response when $\Omega$ is close to the resonant frequency associated with the eigenmode, as seen in \Cref{fig:ampli_delay}. Different magnitudes of force $F$ are shown. When the force is smaller, the response is much greater, thereby allowing the model to capture a very large range of forcing amplitudes with only relatively small variations in acoustic pressure. 

We can also observe (by solving for $\psi_k$) that a phase delay of half a cycle is accumulated as we cross the resonant frequency. The group delay, the time required for information to be delivered, is then given by the derivative $-\mathrm{d}\psi/\mathrm{d}\Omega$. It is observed that delays of several milliseconds are observed in the vicinity of resonance (\Cref{fig:ampli_delay}).

\begin{figure}
	\centering
	\begin{minipage}{.47\textwidth}
		\centering
		\includegraphics[width=0.9\linewidth]{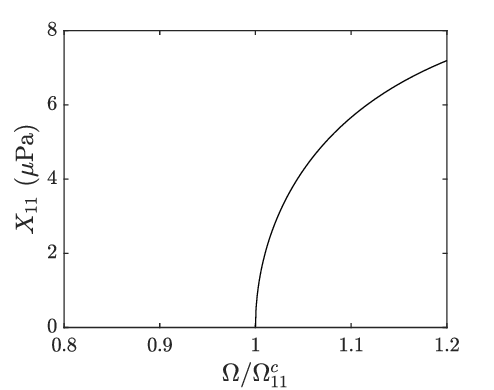}
		\captionof{figure}{The backbone curve of the single-mode equation. Thanks to the non-linearity, the natural response frequency varies as a function of the amplitude. $\hat{\beta}=10^5\si{\second\per\pascal\squared}$ and $\mu=\mu_{11}^c$.}
		\label{fig:backbone}
	\end{minipage}%
	\hspace{0.5cm}
	\begin{minipage}{.47\textwidth}
		\centering
		\includegraphics[width=0.9\linewidth]{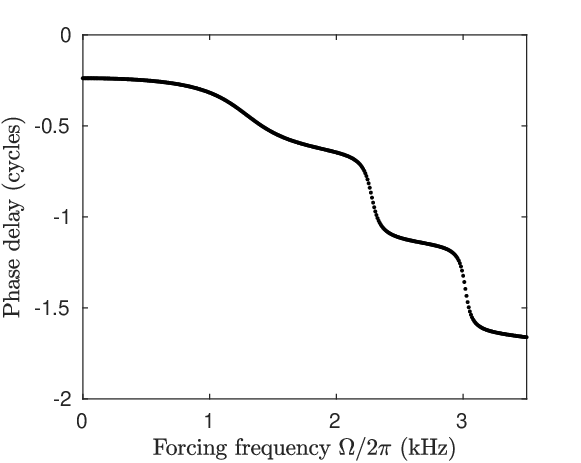}
		\captionof{figure}{The phase delay can, in the fully-coupled system, accumulate to several cycles as the forcing frequency is increased. We study the solution at the centre of the 3\textsuperscript{rd} resonator in response to a sound at 100~dB~SPL.}
		\label{fig:phase_delay}
	\end{minipage}
\end{figure}

When studying relations such as \eqref{eq:ampli_freq}, it becomes apparent that the resonant behaviour occurs slightly away from where is expected (based on the linear system). This is a general property of non-linear systems and can be understood by examining the harmonic response of the unforced non-linear system. Solving \eqref{eq:ampli_freq} in the case $F=0$ gives the relationship shown in \Cref{fig:backbone}, known as a \emph{backbone curve} \cite{wagg2016nonlinear}. We see that with increasing amplitude the natural harmonic response of the non-linear system is perturbed away from the resonant frequency.

\subsection{Fully-coupled system}

Bearing in mind the above analysis of the single-mode approximation \eqref{eq:ode} we now return to the fully-coupled system \eqref{eq:coupled_syst}. Much of the analysis from \Cref{sec:single_mode} can be readily repeated for the matrix system, particularly with the use of numerical schemes for solving non-linear systems of equations, as was used to produce \Cref{fig:full_decomp}. We focus our attention on the elements which tangibly differ from the above discussions.

If we repeat the eigenvalue analysis of \Cref{sec:hopf} we find that a series of Hopf bifurcations take place, at successive parameter values. When we linearise \eqref{eq:coupled_syst} about $\a=0$, since the coupling between modes takes place within the non-linear part of the system, we reach $N$ uncoupled linear systems each of which has Jacobian of the form \eqref{eq:J}. This means that each time $\mu$ passes one of the critical values $\mu_n^c$, $n=1,\dots,N$ (as defined in \eqref{eq:limit_cycle}) a Hopf bifurcation occurs.

An important feature of the fully-coupled system, which we highlight since it is not the case for the single-mode formulation, is the ability to predict phase delays of more than half a cycle. In uncoupled oscillator systems the phase delay will not exceed half a cycle (\emph{cf.} \Cref{fig:ampli_delay}) however the cochlea is well-known to exhibit delays of several cycles  \cite{bell2012resonance,reichenbach2014physics}. Recalling the decomposition \eqref{eq:soln_decomposition} and the harmonic balance techniques used in \Cref{sec:amplification_selectivity}, the phase delay at a given point $x_0$ is given by the complex argument
\begin{equation} \label{eq:phase_delay}
\arg\left(\sum_{n}R_ne^{i\psi_n}u_n(x_0)\right).
\end{equation}
An example of how \eqref{eq:phase_delay} varies as a function of the harmonic forcing frequency is shown in \Cref{fig:phase_delay}. The alternating cliffs and plateaus are because the delay increases much more quickly in the region of one of the system's resonant frequencies.

\section{Discussion}

We have presented a design for an active acoustic metamaterial that is capable of mimicking the properties of the cochlea. Based on a size-graded array of high-contrast resonators, the structure has similar dimensions to the cochlea and has a resonant spectrum that falls broadly within the range of audible frequencies. This design is able to filter different frequencies in space and, with the introduction of a non-linear amplification term, replicate the fundamental properties of the cochlear amplifier.

A modal decomposition was used to approach the coupled-resonator problem. A thorough analysis of a reduced version of the system was undertaken in order to demonstrate its key properties. The aspects which differ most significantly in the fully-coupled case, such as the accumulation of large phase delays, were accounted for.

A final, interesting, aside is that, while the array of compressible resonators considered in this work is presented as a model for the cochlea on the grounds that it can simulate the oscillations of the basilar membrane, there is evidence that the cochlea itself contains compressible elements \cite{shera1992empirical}. The details of this compressibility, and its relevance to cochlear function, are not yet clear.

The code developed for this study is available online at \newline \url{https://github.com/davies-b/hopf_active_cochlea}.

\section*{Acknowledgement}
The authors would like to thank Andrew Bell and A. James Hudspeth  for insightful comments made on an early version of this manuscript.



\appendix
\section{Multipole expansion method} \label{app:multipole}

The layer-potential representation \eqref{eq:layer_potential_representation} 
%
reduces the Helmholtz problem \eqref{eq:helmholtz_equation} to finding density functions, $\phi$ and $\psi$, such that the two transmission conditions on $\D$ are satisfied. That is, we wish to find $\phi,\psi\in L^2(\D)$ such that
\begin{equation} \label{eq:matrix_eqn}
\begin{pmatrix}
\S_D^{k_b} & -\S_D^{k} \\
\partial_\nu^-[\S_D^{k_b}] & -\delta \partial_\nu^+ [\S_D^{k}]
\end{pmatrix}
\begin{pmatrix}
\phi \\ \psi
\end{pmatrix}
=
\begin{pmatrix}
0 \\ 0
\end{pmatrix},
\end{equation}
where equality holds as elements of $L^2(\D)$ and the notation $\partial_\nu^-$ and $\partial_\nu^+$ denotes the derivative in the direction of the normal to $\D$, from the inside and outside of $D$, respectively. More details on the use of layer potentials in solving scattering problems can be found in \emph{e.g.} \cite{ammari2009layer, ammari2018mathematical}.

Since we are interested in the case of circular resonators, $\phi$ and $\psi$ are, on each $\D_n$, $2\pi$-periodic functions of $\theta_n$ where $(r_n,\theta_n)$ denotes a polar coordinate system about on the centre of $D_n$. Such functions admit Fourier expansions of the form
\begin{equation*}
\phi|_{\D_n}=\sum_{m\in\mathbb{Z}} a_m^n e^{im\theta_n},
\end{equation*}
for coefficients $a_m^n$, and similarly for $\psi$. The reason such an expansion is useful is that $\S_{D_n}^k[e^{im\theta_n}]$ has an explicit representation, shown in \cite{ammari2017subwavelength} to be given by
\begin{equation} \label{eq:multipole_formula}
\S_{D_n}^k[e^{im\theta_n}] = \begin{cases}
c_n J_m(k R_n)H_m^{(1)}(kr_n)e^{im\theta_n} & r_n>R_n, \\
c_n H_m^{(1)}(k R_n)J_m(kr_n)e^{im\theta_n} & r_n\leq R_n,
\end{cases}
\end{equation}
where $J_m$ and $H_m^{(1)}$ are the Bessel and Hankel functions of the first kind, respectively, $c_n=-\frac{i\pi R_n}{2}$ and $R_n$ is the radius of $D_n$.

In order to apply this method to the case of $N\in\mathbb{N}$ resonators we also require an expression for $\S_{D_{n'}}^k[e^{im\theta_{n}}]$, where $n\neq n'$. This is achieved through the use of Graf's addition formula \cite{ammari2018mathematical}, which says that for any $x,y\in\mathbb{R}^2$ such that $|x|>|y|$, it holds that the Helmholtz green's function $\Gamma^k$ is given by
\begin{equation} \label{eq:graf}
\Gamma^k(x-y) = -\frac{i}{4} \sum_{l\in\mathbb{Z}} H_l^{(1)}(k|x|) J_l(k|y|) e^{il(\vartheta_x-\vartheta_y)},
\end{equation}
where $x=(|x|,\vartheta_x)$ and $y=(|y|,\vartheta_y)$ are polar representations around a common origin.

Finally, we make the identification $L^2(\D) \cong L^2(\D_1) \times \ldots\times L^2(\D_N)$ and decompose the single layer potential as
\begin{equation} \label{eq:S_matrix}
\S_D^{k} = \begin{pmatrix} 
\S_{D_1}^{k} & \S_{D_2}^k\big|_{\D_1} & \dots & \S_{D_N}^k\big|_{\D_1}  \\[1em]  
\S_{D_1}^k\big|_{\D_2} & \S_{D_2}^{k} & \dots & \S_{D_N}^k\big|_{\D_2}\\
\vdots & \vdots & \ddots & \vdots \\
\S_{D_1}^k\big|_{\D_N} & \S_{D_2}^k\big|_{\D_N} & \dots & \S_{D_N}^{k}
\end{pmatrix},
\end{equation}
where $\S_{D_{n'}}^k\big|_{\D_n}: L^2(\D_{n'}) \rightarrow L^2(\D_n)$ is the evaluation of $\S_{D_{n'}}^{k}$ on $\D_n$. Let $z_{n}^{n'}$ be the vector from the centre of $D_n$ to that of $D_{n'}$, then the off-diagonal terms in \eqref{eq:S_matrix} take the form $\S_{D_{n'}}^k[e^{im\theta_n}]=\S_{z_{n}^{n'}+D_{n}}^k[e^{im\theta_n}]$. The addition of $z_{n}^{n'}$ within the integrand can then be decomposed using \eqref{eq:graf}.

The derivatives appearing in \eqref{eq:matrix_eqn} can be handled similarly, based on the expressions
\begin{align*}
\partial_\nu^-\big[\S_{D_n}^k[e^{im\theta_n}]\big] &= c_n k J_m'(k R_n)H_m^{(1)}(kR_n)e^{im\theta_n} \\
\partial_\nu^+\big[\S_{D_n}^k[e^{im\theta_n}]\big] &= c_n k J_m(k R_n){H_m^{(1)}}'(kR_n)e^{im\theta_n},
\end{align*}
which can be derived by differentiating \eqref{eq:multipole_formula}.

Finally, we truncate the Fourier basis on each $\D_n$, using only $\{e^{im\theta_{n}}:n=-M,\dots,M\}$ for some $M\geq0$, to reach an approximate matrix representation for \eqref{eq:matrix_eqn}.

\end{document}